\documentclass{amsart}

\usepackage[leqno]{amsmath}
\usepackage{latexsym,amsthm,amsfonts,amssymb,hyperref}
\usepackage[all]{xy} \SelectTips{eu}{}

\newcommand{\numberseries}{\mdseries}   
\newlength{\thmtopspace}                
\newlength{\thmbotspace}                
\newlength{\thmheadspace}               
\newlength{\thmindent}                  
\newtheoremstyle{bfupright head,upright body}
                {\thmtopspace}{\thmbotspace}
                {\upshape}{\thmindent}{\bfseries}{.}{\thmheadspace}
                {{\numberseries \thmnumber{(#2) }}\thmnote{#3}}

\newtheoremstyle{fixed bf head,slanted body}
                {\thmtopspace}{\thmbotspace}{\slshape}
                {\thmindent}{\bfseries}{.}{\thmheadspace}
                {{\numberseries \thmnumber{(#2) }}\thmname{#1}\thmnote{ (#3)}}

\newtheoremstyle{fixed bf head,upright body}
                {\thmtopspace}{\thmbotspace}{\upshape}
                {\thmindent}{\bfseries}{.}{\thmheadspace}
                {{\numberseries \thmnumber{(#2) }}\thmname{#1}\thmnote{ (#3)}}

\newtheoremstyle{numbered paragraph}
                {\thmtopspace}{\thmbotspace}{\upshape}
                {\thmindent}{\upshape}{}{0pt}
                {{\numberseries \thmnumber{(#2) }}}

\theoremstyle{bfupright head,upright body}
\newtheorem{res}{}[section]             \newtheorem*{res*}{}
               \newtheorem*{bfhpg*}{}

\theoremstyle{fixed bf head,slanted body}
\newtheorem{thm}[res]{Theorem}          \newtheorem*{thm*}{Theorem}
\newtheorem{prp}[res]{Proposition}      \newtheorem*{prp*}{Proposition}
\newtheorem{cor}[res]{Corollary}        \newtheorem*{cor*}{Corollary}
\newtheorem{lem}[res]{Lemma}            \newtheorem*{lem*}{Lemma}

\theoremstyle{fixed bf head,upright body}
       \newtheorem*{dfn*}{Definition}
\newtheorem{exa}[res]{Example}       \newtheorem*{exa*}{Example}
           \newtheorem*{rmk*}{Remark}

\theoremstyle{numbered paragraph}
\newtheorem{ipg}[res]{}

\newlength{\thmlistleft}        
\newlength{\thmlistright}       
\newlength{\thmlistpartopsep}   
\newlength{\thmlisttopsep}      
\newlength{\thmlistparsep}      
\newlength{\thmlistitemsep}     

\newcounter{prt}
\newenvironment{prt}{\begin{list}{\upshape (\alph{prt})}%
    {\usecounter{prt}%
      \setlength{\leftmargin}{\thmlistleft}%
      \setlength{\labelwidth}{\thmlistleft}%
      \setlength{\rightmargin}{\thmlistright}%
      \setlength{\partopsep}{\thmlistpartopsep}%
      \setlength{\topsep}{\thmlisttopsep}%
      \setlength{\parsep}{\thmlistparsep}%
      \setlength{\itemsep}{\thmlistitemsep}}}%
  {\end{list}}%

\newenvironment{prf*}[1][Proof]{%
  \begin{proof}[\bf #1]
    \setcounter{equation}{0}
    }
  {\end{proof}
}
\renewcommand{\eqref}[1]{\pgref{eq:#1}}
\newcommand{\pgref}[1]{(\ref{#1})}
\newcommand{\thmref}[2][Theorem~]{#1\pgref{thm:#2}}
\newcommand{\prpref}[2][Proposition~]{#1\pgref{prp:#2}}
\newcommand{\corref}[2][Corollary~]{#1\pgref{cor:#2}}
\newcommand{\lemref}[2][Lemma~]{#1\pgref{lem:#2}}
\newcommand{\thmcite}[2][?]{\cite[Theorem~#1]{#2}}

\newcommand{\prpcite}[2][?]{\cite[Proposition~#1]{#2}}
\newcommand{\corcite}[2][?]{\cite[Corollary~#1]{#2}}
\newcommand{\seccite}[2][?]{\cite[Section~#1]{#2}}

\def\urltilda{\kern -.15em\lower .7ex\hbox{\~{}}\kern .04em} 
\newcommand{\set}[2][\mspace{1mu}]{\{#1 #2 #1\}}

\newcommand{\kk}{\Bbbk}
\newcommand{\qtext}[1]{\quad\text{#1}\quad}

\newcommand{\qand}{\qtext{and}}

\newcommand{\deq}{\:=\:}
\newcommand{\dle}{\:\le\:}

\newcommand{\is}{\cong}

\newcommand{\lla}{\longleftarrow}

\newcommand{\xla}[2][]{\xleftarrow[#1]{\;#2\;}}

\renewcommand{\H}[2][]{\operatorname{H}_{#1}(#2)}

\newcommand{\rnk}[2][k]{\operatorname{rank}_{#1}#2}

\newcommand{\Tor}[4][Q]{\operatorname{Tor}^{#1}_{#2}(#3,#4)}

\newcommand{\fn}{\mathfrak{n}}
\newcommand{\ff}{\mathfrak{f}}

\newcommand{\calI}{\mathcal{I}}
\newcommand{\calJ}{\mathcal{J}}
\newcommand{\calN}{\mathcal{N}}
\newcommand{\calQ}{\mathcal{Q}}

\numberwithin{equation}{res}

\renewcommand{\le}{\leqslant}
\renewcommand{\ge}{\geqslant}

\begin{document}

\title[Free resolutions of Dynkin format]{Free resolutions of Dynkin
  format\\ and the licci property of grade $3$ perfect ideals}

\author[L.\,W. Christensen]{Lars Winther Christensen}

\address{Texas Tech University, Lubbock, TX 79409, U.S.A.}

\email{lars.w.christensen@ttu.edu}

\urladdr{http://www.math.ttu.edu/\urltilda lchriste}

\author[O. Veliche]{Oana Veliche}

\address{Northeastern University, Boston, MA~02115, U.S.A.}

\urladdr{https://web.northeastern.edu/oveliche}

\email{o.veliche@northeastern.edu}

\author[J. Weyman]{Jerzy Weyman}

\address{University of Connecticut, Storrs, CT~06269, U.S.A.}

\email{jerzy.weyman@uconn.edu}

\urladdr{http://www.math.uconn.edu/\urltilda weyman}

\thanks{This work is part of a body of research that started during
  the authors' visit to MSRI in Spring 2013 and continued during a
  months-long visit by L.W.C.\ to Northeastern University; the
  hospitality of both institutions is acknowledged with
  gratitude. L.W.C.\ was partly supported by NSA grant H98230-14-0140
  and Simons Foundation collaboration grant 428308, and J.W.\ was
  partly supported by NSF DMS grant 1400740.}

\date{17 January 2019}

\keywords{Free resolution, licci ideal, linkage, perfect ideal}

\subjclass[2010]{13C40; 13D02.}

\begin{abstract}
  Recent work on generic free resolutions of length 3 attaches to
  every resolution a graph and suggests that resolutions whose
  associated graph is a Dynkin diagram are distinguished.  We
  conjecture that in a regular local ring, every grade $3$ perfect
  ideal whose minimal free resolution is distinguished in this way is
  in the linkage class of a complete intersection.
\end{abstract}

\maketitle

\thispagestyle{empty}

\section{Introduction}
\label{sec:introduction}

\noindent
Let $Q$ be a commutative Noetherian ring. Quotient rings of $Q$ that
have projective dimension at most $3$ as $Q$-modules have been and
still are investigated vigorously. A recent development is Weyman's
\cite{JWm} construction of generic rings for resolutions of length
$3$. For a free resolution,
\begin{equation*}
  F_0 \xla{\partial_1} F_1 \xla{\partial_2} F_2 \xla{\partial_3} F_3 \lla 0\:,
\end{equation*}
the \emph{format} is the quadruple $(f_0,f_1,f_2,f_3)$ with
$f_i=\rnk[]{F_i}$. To each resolution format Weyman associates a graph
and a generic ring, and he proves that the generic ring is Noetherian
exactly if the graph is a Dynkin diagram.  This suggests that free
resolutions of formats corresponding to Dynkin diagrams play a special
role; here we explore what that role might be in the context of
linkage.

Let $Q$ be local and $I$ in $Q$ be a perfect ideal of grade $3$. Let
$m$ be the minimal number of generators of $I$, the minimal free
resolution of $Q/I$ over $Q$ then has format $\ff = (1,m,m+n-1,n)$ for
$m\ge 3$ and some integer $n \ge 1$ which, when $Q$ is regular, is
referred to as the \emph{type} of the ring $Q/I$. For convenience we
refer to the tuple $\ff$ as the \emph{resolution format}, or simply
the \emph{format,} of $I$. The graph associated to $\ff$ is
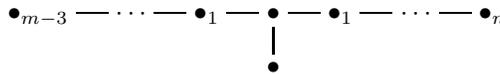
\begin{equation*}
  \xymatrix@C=1pc@R=.75pc{
    \bullet_{m-3} \ar@{-}[r] & 
    \cdots \ar@{-}[r] & 
    \bullet_1 \ar@{-}[r] & \bullet \ar@{-}[r] & 
    \bullet_1 \ar@{-}[r] &
    \cdots \ar@{-}[r] & \bullet_{n} \\
    &&&\bullet \ar@{-}[u]
  }
\end{equation*}
A complete intersection ideal has resolution format $(1,3,3,1)$, and
the corresponding graph is the Dynkin diagram $\mathrm{A}_3$. A
Gorenstein ideal that is not complete intersection has format
$(1,m,m,1)$ for some odd $m\ge 5$, and the associated graph is the
Dynkin diagram $\mathrm{D}_m$. There is another way to obtain Dynkin
diagrams of type $\mathrm{D}$, namely from almost complete
intersection ideals, which have formats $(1,4,n+3,n)$, corresponding
to $\mathrm{D}_{n+3}$, for $n\ge 2$. It is elementary to verify that
$\mathrm{E}_6$, $\mathrm{E}_7$, and $\mathrm{E}_8$ are the only other
Dynkin diagrams that can possibly arise in this manner. In the
remainder of this paper, we say that a (resolution) format is
\emph{Dynkin} if the corresponding graph is a Dynkin diagram.

If $Q$ is regular, then every grade $2$ perfect ideal $I$ in $Q$ is
\emph{licci}, that is, in the linkage class of a complete intersection
ideal. Not every perfect ideal of grade $3$ is licci, and in a certain
sense licci ideals of grade $4$ are rare; see Miller and
Ulrich~\cite{MMlBUl86}. Works of Kunz~\cite{EKz74} and
Watanabe~\cite{JWt73} show that a grade $3$ perfect ideal is licci if
it is Gorenstein or almost complete intersection; these results
predate the term licci, but Buchsbaum and Eisenbud interpret them in
the introduction to \cite{DABDEs77}. The purpose of this paper is to
motivate the following:

\begin{bfhpg*}[Conjecture]
  Let $Q$ be a regular local ring and $\ff = (1,m,m+n-1,n)$ be a
  resolution format realized by some grade $3$ perfect ideal in $Q$.
  \begin{prt}
  \item[\textbf{I}\,] If $\ff$ is not Dynkin, then there exists a
    grade $3$ perfect ideal in $Q$ of format $\ff$ that is not licci.
  \item[\textbf{II}\,] If $\ff$ is Dynkin, then every grade $3$
    perfect ideal in $Q$ of format $\ff$ is licci.
  \end{prt}
\end{bfhpg*}
\noindent
The paper is organized as follows. In Section \ref{sec:linkage} we
recall basics on linkage and prove a few technical results. In Section
\ref{sec:I} we collect evidence for part I of the conjecture; in fact,
we prove it for local rings that are obtained by localizing a
polynomial algebra at the irrelevant maximal ideal. Some of our
evidence for part II was already discussed above, and further evidence
is provided in Section \ref{sec:II}.

\section{Linkage of grade $3$ perfect ideals}
\label{sec:linkage}

\noindent
To goal of this section, other than to recall the language of linkage,
is to establish a proposition that allows some control over the
resolution formats of directly linked ideals. To this end we need a
folklore application of Prime Avoidance for which we did not find a
citable reference.\footnote{Note added in proof: A, not widely
  available, reference is Kaplansky \thmcite[124]{Kap}.}

\begin{lem}
  \label{lem:prav}
  Let $Q$ be a commutative ring and $P_1,\ldots,P_s$ be prime ideals
  in $Q$. For every ideal $J$ and element $x$ in $Q$ with
  $J + (x) \not\subseteq P_1\cup P_2\cup\ldots\cup P_s$ there exists
  an element $v\in J$ with $v+x \notin P_1\cup P_2\cup\ldots\cup P_s.$
\end{lem}

\begin{prf*} 
  Without loss of generality we may assume that there are no
  containments among the prime ideals $P_1,\ldots,P_s$ and that they
  are ordered such that there is an $r$ between $0$ and $s$ with
  $x\in P_1,\ldots,P_r$ and $x\notin P_{r+1},\ldots ,P_s$. Notice that
  if $r=0$, then one can take $v=0$.

  Let $r>0$. From the assumption on $J + (x)$ one gets
  $J\not\subseteq P_1\cup\ldots\cup P_r$; choose an element
  $u\in J\setminus (P_1\cup\ldots\cup P_r)$. If $r=s$, then one can
  take $v=u$. In the case $r<s$ one has
  $P_{r+1}\cap\ldots \cap P_s\not\subseteq P_1\cup\ldots\cup P_r$ as,
  indeed, by Prime Avoidance containment would imply
  $P_{r+1}\cap\ldots \cap P_s \subseteq P_i$ for some $i$ and then one
  would have $P_j\subseteq P_i$ for some $j\ne i$, which contradicts
  the assumption on $P_1,\ldots,P_s$.  Now choose an element
  $y\in P_{r+1}\cap\ldots \cap P_s \setminus P_1\cup\ldots\cup
  P_r$. It is elementary to verify that $uy+x$ is not in
  $P_1\cup P_2\cup\ldots\cup P_s$.
\end{prf*}

In the rest of this section, we keep the setup and notation close to
\cite{AKM-88} by Avramov, Kustin, and Miller; the reader may see
\seccite[1]{AKM-88} for more background.

Let $I$ be a perfect ideal of grade $3$ in a local ring $Q$. An ideal
$J \subseteq Q$ is said to be \emph{directly linked} to $I$ if there
exists a regular sequence $x_1, x_2, x_3$ of elements in $I$ with
$J = (x_1, x_2, x_3):I$. The ideal $J$ is then also a perfect ideal of
grade $3$.  An ideal $J$ is said to be \emph{linked} to $I$ if there
exists a sequence of ideals $I=J_0, J_1,\ldots,J_n = J$ such that
$J_{i+1}$ is directly linked to $J_i$ for each $i=0,\ldots ,n-1$.

\begin{prp}
  \label{prp:link}
  Let $Q$ be a local ring and $I \subseteq Q$ a grade $3$ perfect
  ideal of resolution format $(1,m,m+n-1,n)$.
  \begin{prt}
  \item If $m \ge 4$ holds, then there exists a grade $3$ perfect
    ideal of resolution format $(1,n+3,m+n,m-2)$ that is directly
    linked to $I$.
  \item If $m \ge 5$ holds, then there exists a grade $3$ perfect
    ideal of resolution format $(1,n+3,m+n-1,m-3)$ that is directly
    linked to $I$.
  \end{prt}
\end{prp}

\begin{prf*}
  Denote by $\fn$ the maximal ideal of $Q$ and by $k$ the residue
  field $Q/\fn$. Let $F_\bullet$ be the minimal free resolution of
  $Q/I$. For a regular sequence $x_1,x_2,x_3$ in $I$ the corresponding
  Koszul complex is denoted by $K_\bullet$. The canonical surjection
  $Q/(x_1,x_2,x_3) \to Q/I$ lifts to a morphism
  $\varphi_\bullet \colon K_\bullet \to F_\bullet$. As recalled in
  \seccite[1]{AKM-88}, the dual of the mapping cone of
  $\varphi_\bullet$ yields a, not necessarily minimal, free resolution
  of $Q/J$. Further, the ranks of the modules in the minimal free
  resolution $G_\bullet$ of $Q/J$ are given by
  \begin{align*}
    \rnk[Q]{G_1} &= n + 3 - \rnk{(\varphi_3\otimes k)} - \rnk{(\varphi_2\otimes k)}\\
    \rnk[Q]{G_2} &= m + n + 2 - \rnk{(\varphi_2\otimes k)} - \rnk{(\varphi_1\otimes k)}\\
    \rnk[Q]{G_3} &= m - \rnk{(\varphi_1\otimes k)}\:.
  \end{align*}
  Set $A_\bullet = \Tor{\bullet}{Q/I}{k}$.  The morphism
  $\varphi_\bullet \otimes k\colon \H{K_\bullet} \to A_\bullet$ is a
  morphism of graded-commutative $k$-algebras.

  (a): By assumption, $I$ is not generated by a regular sequence, so
  by \thmcite[2.1]{AKM-88} one has $(A_1)^3=0$ and, therefore,
  $\rnk{(\varphi_3\otimes k)}=0$.  It now suffices to show that one
  can choose a regular sequence $x_1,x_2,x_3$ such that
  $\rnk{(\varphi_1\otimes k)}=2$ and $\rnk{(\varphi_2\otimes k)}=0$
  hold. By \thmcite[2.1]{AKM-88} one can choose minimal generators
  $x'_1$ and $x'_2$ of $I$ such that the corresponding cycles $z'_1$
  and $z'_2$ in $F_1$ satisfy $[z'_1]\cdot[z'_2] = 0$ in
  $A_\bullet$. Let $P_1,\ldots,P_s$ be the associated prime ideals of
  $Q$. As $I$ has grade $3$, the ideal $\fn I + (x'_1)$ is not
  contained in the union $\bigcup_{i=1}^s P_i$. It now follows from
  \lemref{prav} that there exists an element $v_1$ in $\fn I$ such
  that $x_1 = v_1 + x_1'$ is not in $\bigcup_{i=1}^s P_i$, i.e.\ not a
  zero-divisor. Similarly, there exists a $v_2\in\fn I$ such that
  $x_2 = v_2 + x_2'$ is not in the union of the associated primes of
  the $Q$-module $Q/(x_1)$. Thus, $x_1,x_2$ form a regular sequence of
  minimal generators of $I$, and the corresponding cycles $z_1$ and
  $z_2$ in $F_1$ satisfy $[z_1]=[z'_1]$ and $[z_2]=[z'_2]$ in
  $A_1$. Finally there exists, again by \lemref{prav}, an element
  $x_3$ in $\fn I$ and not in the union of the associated primes of
  the $Q$-module $Q/(x_1,x_2)$. With this regular sequence
  $x_1,x_2,x_3$ one has $\rnk{(\varphi_1\otimes k)}=2$, as the cycle
  $z_3$ in $F_1$ corresponding to $x_3$ satisfies $[z_3]=0$ in
  $A_1$. In particular, one has
  $[z_1]\cdot [z_3] = 0 = [z_2]\cdot [z_3]$ and hence
  $\rnk{(\varphi_2\otimes k)}=0$.

  (b): By assumption, $I$ is not generated by a regular sequence, so
  by \thmcite[2.1]{AKM-88} one has $(A_1)^3=0$ and, therefore,
  $\rnk{(\varphi_3\otimes k)}=0$. It now suffices to show that one can
  choose a regular sequence $x_1,x_2,x_3$ such that
  $\rnk{(\varphi_1\otimes k)}=3$ and $\rnk{(\varphi_2\otimes k)}=0$
  hold. By \thmcite[2.1]{AKM-88} one can choose minimal generators
  $x'_1$, $x'_2$, and $x'_3$ of $I$ such that the corresponding cycles
  $z'_1, z'_2, z'_3 \in F_1$ satisfy $[z'_i]\cdot[z'_j] = 0$ in
  $A_\bullet$. (Notice that the assumption $m \ge 5$ is needed in
  order to choose these elements in case the multiplicative structure
  is TE; for the other possible structures, B, G, and H, such elements
  can be chosen as long as $m$ is at least $4$.) As in the proof of
  part (a) there exist elements $v_1,v_2,v_3\in \fn I$ such that
  $v_1 + x'_1,v_2 + x'_2,v_3 + x'_3$ is a regular sequence in $I$ with
  the desired properties.
\end{prf*}

\begin{cor}
  \label{cor:link}
  Let $Q$ be a local ring and $I \subseteq Q$ a grade $3$ perfect
  ideal of resolution format $(1,m,m+n-1,n)$.
  \begin{prt}
  \item If $m \ge 4$ holds, then there exists a grade $3$ perfect
    ideal of resolution format $(1,m+1,m+n+1,n+1)$ that is linked to
    $I$.
  \item If $m \ge 5$ holds, then there exists a grade $3$ perfect
    ideal of resolution format $(1,m,m+n,n+1)$ that is linked to $I$.
  \item If $m\ge 4$ and $n \ge 2$ holds, then there exists a grade $3$
    perfect ideal of resolution format $(1,m+1,m+n,n)$ that is linked
    to $I$.
  \end{prt}
\end{cor}

\begin{prf*}
  Part (a) follows from two consecutive applications of
  \prpref{link}(a). Part (b) follows from an application of
  \prpref[]{link}(b) followed an application of
  \prpref[]{link}(a). Finally, part (c) follows by application of
  \prpref[]{link}(a) and \prpref[]{link}(b) in that order.
\end{prf*}

Let $I \subseteq Q$ be a perfect ideal of grade $3$. For every regular
sequence $x_1,x_2,x_3$ in $I$ one has
$(x_1,x_2,x_3):((x_1,x_2,x_3):I)=I$; see Golod \cite{ESG80}. Thus if
an ideal $J$ is directly linked to $I$, then $I$ is also directly
linked to $J$. It follows that ``being linked'' is an equivalence
relation. The ideal $I$ is called \emph{licci} if it is in the linkage
class of a complete intersection ideal.

\section{Evidence for Part I}
\label{sec:I}

\noindent
By a Dynkin diagram we always mean a simply laced Dynkin diagram,
i.e.\ a diagram from the ADE classification. Recall from the
Introduction that
\begin{equation*}
  \xymatrix@C=1pc@R=.75pc {
    \bullet_{m-3} \ar@{-}[r] & 
    \cdots \ar@{-}[r] & 
    \bullet_1 \ar@{-}[r] & \bullet \ar@{-}[r] & 
    \bullet_1 \ar@{-}[r] &
    \cdots \ar@{-}[r] & \bullet_{n} \\
    &&&\bullet \ar@{-}[u]
  }
\end{equation*}
is the graph associated to the resolution format $(1,m,m+n-1,n)$. The
Dynkin formats that can be realized by grade $3$ perfect ideals are
\begin{equation*}
  \label{Dynkin-formats}
  \tag{3.1}
  \begin{alignedat}{2}
    \mathrm{A}_3 & \quad \text{corresponding to} &\quad & (1, 3, 3, 1)\\
    \mathrm{D}_m & \quad \text{\hspace{2.5em} --- }
    &\quad& (1, m, m, 1) \text{ for odd $m \ge 5$}\\
    \mathrm{D}_{n+3} & \quad \text{\hspace{2.5em} --- }
    &\quad& (1, 4, n+3, n) \text{ for $n \ge 2$}\\
    \mathrm{E}_6 &\quad \text{\hspace{2.5em} --- } &\quad& (1,5,6,2)\\
    \mathrm{E}_7 &\quad \text{\hspace{2.5em} --- } &\quad& (1,6,7,2) \text{ and } (1,5,7,3)\\
    \mathrm{E}_8 &\quad \text{\hspace{2.5em} --- } &\quad& (1,7,8,2)
    \text{ and } (1,5,8,4)
  \end{alignedat}
\end{equation*}
\stepcounter{res} This is straightforward to verify when one recalls
that
\begin{prt}
\item[$\bullet$] $(1,3,3,1)$ is the only possible format with $m=3$;
  this excludes the formats $(1,3,n+2,n)$ corresponding to
  $\mathrm{A}_{n+2}$ for $n\ge 2$.
\item[$\bullet$] The minimal number of generators of grade $3$
  Gorenstein ideal is odd, see \cite{DABDEs77}. This excludes the
  formats $(1,m,m,1)$ corresponding to $\mathrm{D}_m$ for even
  $m\ge 4$.
\end{prt}

\begin{thm}
  \label{thm:notdynkin}
  Let $\kk$ a field and $e\ge 3$; set
  $Q = \kk[X_1,\ldots,X_e]_{(X_1,\ldots,X_e)}$. For every resolution
  format $\ff = (1,m,m+n-1,n)$ with $m\ge 3$ and $n\ge 1$ that is not
  Dynkin there exists a grade $3$ perfect ideal that has resolution
  format $\ff$ and is not licci.
\end{thm}

It is implicit in this statement that every format $(1,m,m+n-1,n)$
that is not Dynkin is realized by a degree $3$ perfect ideal in the
ring $Q$. A proof was sketched in the last section of \cite{JWm}, we
provide the details in the Appendix.

The proof of \thmref{notdynkin} takes up the bulk of this section;
indeed the Theorem follows from \prpref[Propositions~]{reduction},
\prpref[]{1683}, and \prpref[]{1892}.  It proceeds in two steps: First
we prove that it is enough to exhibit non-licci ideals for two
specific formats, $(1,6,8,3)$ and $(1,8,9,2)$, which are not
Dynkin. Second we produce such examples, in fact entire families of
them.

\begin{prp}
  \label{prp:reduction}
  Let $Q$ be a local ring.  If there is a resolution format $\ff$
  such~that
  \begin{prt}
  \item[$(1)$] there exists a grade $3$ perfect ideal in $Q$ of format
    $\ff$,
  \item[$(2)$] $\ff$ is not Dynkin, and
  \item[$(3)$] every perfect ideal in $Q$ of format $\ff$ is licci
  \end{prt}
  then $(1,6,8,3)$ or $(1,8,9,2)$ is such a format.
\end{prp}

\begin{prf*}
  Assume that there exists at least one resolution format that
  satisfies conditions (1)--(3). Let $\ff = (1,m,m+n-1,n)$ be minimal
  among all such formats with regard to the total rank
  $\beta(\ff) = 2(m+n)$ of the modules in the minimal free
  resolution. Necessarily, one has $m\ge 4$ as the Dynkin format
  $(1,3,3,1)$ is the only format with $m=3$. It follows that the
  format $\ff'=(1,n+2,m+n-2,m-3)$ is Dynkin. Indeed, one has
  $\beta(\ff') = \beta(\ff)-2$ and every perfect ideal of format
  $\ff'$ is by \prpref{link}(a) is directly linked to a perfect ideal
  of format $\ff$ and hence licci.

  Given the Dynkin formats \pgref{Dynkin-formats}, the potential
  minimal formats $\ff$ are as follows:
  \begin{center}
    \begin{tabular}{r|l}
      Dynkin format $\ff'$ & Potential minimal format $\ff$\\ \hline
      $(1, 2l+1, 2l+1, 1)$ for $l \ge 1$ & $(1,4,2l+2,2l-1)$\\
      $(1, 4, l+3, l)$ for $l\ge 2$ & $(1,l+3,l+4,2)$\\
      $(1,5,6,2)$ & $(1,5,7,3)$\\
      $(1,5,7,3)$ & $(1,6,8,3)$\\
      $(1,5,8,4)$ & $(1,7,9,3)$\\    
      $(1,6,7,2)$ & $(1,5,8,4)$\\
      $(1,7,8,2)$ & $(1,5,9,5)$
    \end{tabular}
  \end{center}
  Among the potential minimal formats, $(1,4,2l+2, 2l-1)$ does not
  exist for $l=1$, and for $l \ge 2$ the format is Dynkin
  corresponding to $\mathrm{D}_{2l+2}$. The format
  $\ff_l = (1, l+3, l+4, 2)$ is for $l\in\{2,3,4\}$ Dynkin
  corresponding to $\mathrm{E}_6$, $\mathrm{E}_7$, and
  $\mathrm{E}_8$. Further, the formats $(1,5,7,3)$ and $(1,5,8,4)$ are
  Dynkin corresponding to $\mathrm{E}_7$ and $\mathrm{E}_8$. Thus,
  possible minimal formats that are not Dynkin and such that every
  ideal of the format is licci are $(1,5,9,5)$, $(1,6,8,3)$,
  $(1,7,9,3)$ and $\ff_l$ for $l\ge 5$.

  An ideal of format $(1,5,9,5)$ is by \prpref{link}(b) linked to an
  ideal of format $\ff_5 = (1,8,9,2)$. An ideal of format $(1,6,8,3)$
  is by \corref{link}(c) linked one of format $(1,7,9,3)$. Similarly,
  an ideal of format $\ff_l$ is linked to one of format
  $f_{l+1}$. Thus, if there exists a format that is not Dynkin with
  the property that every ideal of that format is licci, then
  $(1,6,8,3)$ or $(1,8,9,2)$ has that property.
\end{prf*}

To provide examples of ideals of formats $(1,6,8,3)$ and $(1,8,9,2)$
that are not licci, we rely on numerical obstructions found by Huneke
and Ulrich \cite{CHnBUl87}.

\begin{ipg}
  \label{huineq}
  Let $\kk$ be a field and $e\ge 3$; set $\calQ=\kk[X_1,\ldots, X_e]$.
  Let $\calI$ be a homogeneous perfect ideal in $\calQ$ of grade $3$
  with minimal free resolution
  \begin{equation*}
    \calQ \lla  \bigoplus_{i=1}^{m} \calQ(-d_{1,i}) \lla  \bigoplus_{i=1}^{m+n-1} 
    \calQ(-d_{2,i}) \lla \bigoplus_{i=1}^{n} \calQ(-d_{3,i}) \lla 0\:.
  \end{equation*}
  Set $Q = \calQ_{(X_1,\ldots,X_e)}$ and
  $I = \calI_{(X_1,\ldots,X_e)}$. Assuming that one has
  $d_{1,1}\le d_{1,2} \le \cdots \le d_{1,m}$ and
  $d_{3,1}\le \cdots \le d_{3,n}$ it follows from
  \corcite[5.13]{CHnBUl87} that if the inequality
  \begin{equation*}
    d_{3,n} \le 2d_{1,1}
  \end{equation*}
  holds, then the ideal $I$ is not licci.

  Notice that if $Q/I$ is Artinian, i.e.\ $e=3$, then one has
  $Q/I \is \calQ/\calI$ and the inequality says that $3$ plus the
  socle degree of $\calQ/\calI$ does not exceed $2$ times the initial
  degree of $\calI$.
\end{ipg}

\subsection*{Resolution format $(1,6,8,3)$}

\begin{prp}
  \label{prp:1683}
  Adopt the notation from {\rm \pgref{huineq}.} Let $\Phi$ be a
  $2\times 4$ matrix of linear forms in $\calQ$ or a $3\times 3$
  symmetric matrix of linear forms, and let $\calI$ be the ideal of
  $2\times 2$ minors of $\Phi$. The ideal $I$ has resolution format
  $(1,6,8,3)$ and it is not licci.
\end{prp}

\begin{prf*}
  The graded minimal free resolution of $\calQ/\calI$ has the form
  \begin{equation*}
    \calQ \lla \calQ(-2)^6 \lla \calQ(-3)^8 \lla  \calQ(-4)^3 \lla 0\:;
  \end{equation*}
  see Eagon and Nothcott \thmcite[2]{JAEDGN62} and J\'ozefiak
  \thmcite[3.1]{TJz78} or Weyman \prpcite[(6.1.7) and
  thm.~(6.3.1)]{JWm03}.  It follows from \pgref{huineq} that $I$ is
  not licci.
\end{prf*}

The next example illustrates \prpref{1683} and the linkage argument in
the proof of \prpref{reduction}. Recall that an Artinian local ring is
called \emph{level} if the socle is exactly the highest nonvanishing
power of the maximal ideal.

\begin{exa}
  \label{exa:1683}
  Let $\kk$ be a field; set $\calQ = \kk[X,Y,Z]$ and
  $\calN = (X,Y,Z)$.  The ideal $\calN^2$ is generated by the six
  quadratic monomials, and the quotient $\calQ/\calN^2$ is an Artinian
  local level algebra of socle degree 1 and type $3$, so the minimal
  free resolution over $\calQ$ is
  \begin{equation*}
    \calQ \lla \calQ(-2)^6 \lla \calQ(-3)^8 \lla  \calQ(-4)^3 \lla 0\:;
  \end{equation*}
  Notice that $\calN^2$ is the ideal of $2\times 2$ minors of either
  matrix
  \begin{equation*}
    \begin{pmatrix}
      X & Y & Z & 0 \\ 0 & X & Y & Z
    \end{pmatrix}
    \quad\text{or}\quad
    \begin{pmatrix}
      X & Y & Z \\ Y & 0 & X \\ Z & X & Y
    \end{pmatrix}\:.
  \end{equation*}
  The linked ideal
  \begin{equation*}
    \calI = (X^2,Y^2,Z^3):\calN^2 = (X^2,Y^2, XYZ, XZ^2, YZ^2,Z^3)
  \end{equation*}
  defines a local Artinian $\kk$-algebra with graded basis
  \begin{equation*}
    1 \quad x, y, z \quad xy, xz, yz, z^2\:.
  \end{equation*}
  Evidently, there are no linear forms in the socle of $\calQ/\calI$,
  so it is a level algebra of type $4$, whence $\calI$ has resolution
  format $(1,6,9,4)$. Next, the ideal
  \begin{multline*}
    \calJ=(X^3-YZ^2, Y^3-XZ^2,Z^3):\calI \\= (X^3-YZ^2,Y^3-XZ^2,Z^3,
    X^2Y^2,X^2YZ,XY^2Z,XYZ^2)
  \end{multline*}
  defines a local Artinian $\kk$-algebra with graded basis
  \begin{equation*}
    1 \quad x, y, z \quad x^2, xy, xz, y^2, yz, z^2 
    \quad x^2y, x^2z, xy^2, xyz, xz^2,y^2z,yz^2 \quad x^2z^2, y^2z^2\:.
  \end{equation*}
  It is straightforward to verify that the socle of $\calQ/\calJ$ is
  generated by $xyz$, $x^2z^2$, and $y^2z^2$. It follows that $\calJ$
  has resolution format $(1,7,9,3)$.
\end{exa}

\subsection*{Resolution format $(1,8,9,2)$}
Recall that an Artinian $\kk$-algebra is called \emph{compressed} if
it has maximal length among all algebras with the same socle
polynomial; see \pgref{compressed}. This notion has a natural
generalization to Cohen--Macaulay algebras, see Fr\"oberg and Laksov
\cite{RFrDLk84}, and that is the one we use in the proposition below.

\begin{prp}
  \label{prp:1892}
  Adopt the notation from {\rm \pgref{huineq}.} Assume that the graded
  $\kk$-algebra $\calQ/\calI$ is compressed with with socle polynomial
  $2t^3$; that is,
  \begin{equation*}
    \dim_\kk (\Tor[\calQ]{3}{\calQ/\calI}{\kk})_i = 
    \begin{cases}
      2 & \text{for $i=6$}\\
      0 & \text{for $i\ne 6\:.$}
    \end{cases}
  \end{equation*}
  The ideal $I$ has resolution format $(1,8,9,2)$ and it is not licci.
\end{prp}

\begin{prf*}
  By \prpcite[6]{RFrDLk84} the Hilbert series of $\calQ/\calI$ is
  \begin{equation*}
    \frac{1 + 3t + 6t^2 + 2t^3}{(1-t)^{e-3}}\:.
  \end{equation*}
  Further, it follows from the equality $\binom{3}{2} = 2\binom{3}{1}$
  that $\calQ/\calI$ is extremely compressed, see
  \cite[p.~133]{RFrDLk84}. Now it follows from \prpcite[16]{RFrDLk84}
  that the graded free resolution of $\calQ/\calI$ has the form
  \begin{equation*}
    \calQ \lla \calQ(-3)^8 \lla \calQ(-4)^9 \lla  \calQ(-6)^2 \lla 0\:;
  \end{equation*}
  The conclusion now follows from \pgref{huineq}.
\end{prf*}

By work of Boij and Laksov \thmcite[(3.4)]{MBjDLk94}, \prpref{1892}
applies to generic Artinian level algebras of embedding dimension $3$,
socle degree $3$, and type $2$.

The next example illustrates \prpref{1892} and the linkage argument in
the proof of \prpref{reduction}.

\begin{exa}
  \label{exa:1892}
  Let $\kk$ be a field. The ideal
  \begin{equation*}
    \calI =(X^3, X^2Y+YZ^2, X^2Z+XYZ, XY^2+XYZ,XZ^2,Y^3,Y^2Z,Z^3)
  \end{equation*}
  in $\kk[X,Y,Z]$ defines an Artinian local $\kk$-algebra with graded
  basis
  \begin{equation*}
    1 \quad x, y, z \quad x^2, xy, xz, y^2, yz, z^2 \quad xyz, yz^2\:.
  \end{equation*}
  It is elementary to check that there are no quadratic forms in the
  socle of $\calQ/\calI$, so it is a compressed level algebra of type
  $2$, and it follows that $\calI$ has format $(1,8,9,2)$.

  The linked ideal
  \begin{equation*}
    \calJ = (X^3,Y^3,Z^3):\calI = (X^3, X^2Y-YZ^2, XYZ - XZ^2-Y^2Z,Y^3,Z^3)
  \end{equation*}
  defines an Artinian local $\kk$-algebra with graded basis
  \begin{equation*}
    1 \quad x, y, z \quad x^2, xy, xz, y^2, yz, z^2 \quad x^2z, xy^2,xz^2,y^2z,yz^2\:.
  \end{equation*}
  Again it is elementary to verify that there are no quadratic forms
  in the socle, so it is a level algebra of type $5$; in particular,
  the resolution format of $\calJ$ is $(1,5,9,5)$.

  Finally, notice that under the flat extension
  $\kk[X,Y,Z] \to \kk[X,Y,Z,X_4,\ldots,X_e]$ the ideal $\calI$ extends
  to an ideal that satisfies the hypothesis of \prpref{1892}.
\end{exa}

\section{Evidence for Part II}
\label{sec:II}

\noindent
As discussed in the Introduction, it is known that grade $3$ perfect
ideals of resolution format $(1,4,n+3,n)$ for $n\ge 2$---almost
complete intersection ideals---and ideals of format $(1,m,m,1)$ for
odd $m\ge 3$---Gorenstein ideals---are licci.  The remaining Dynkin
formats are $(1,5,6,2)$ corresponding to $\mathrm{E}_6$, $(1,6,7,2)$
and $(1,5,7,3)$ corresponding to $\mathrm{E}_7$, and $(1,7,8,2)$ and
$(1,5,8,4)$ corresponding to $\mathrm{E}_8$; see
\pgref{Dynkin-formats}.

In the previous section, we proved the existence of ideals that are
not licci by invoking the numerical obstructions in \pgref{huineq}. In
this section, we prove that homogeneous grade $3$ perfect ideals in
$\kk[X_1,X_2,X_3]$ of resolution formats corresponding to
$\mathrm{E}_6$, $\mathrm{E}_7$, or $\mathrm{E}_8$ avoid this
obstruction.

\begin{ipg}
  \label{compressed}
  Let $\calQ$ be a polynomial algebra of embedding dimension $e$. For
  a finitely generated graded $\calQ$-module $M$, let $h_M(\cdot)$
  denote the Hilbert function. Let $\calI$ be a homogeneous ideal in
  $\calQ$ such that $\calQ/\calI$ is Artinian with socle degree $s$;
  its socle polynomial $\sum_{j=1}^sc_jt^j$ is the Hilbert series of
  the socle of $\calQ/\calI$. For every $i\ge 0$ there is an
  inequality
  \begin{align*}
    h_{\calQ/\calI}(i) &\dle \min\set{h_\calQ(i),\sum_{j=i}^sc_jh_\calQ(j-i)} \\
                       &\deq \min\left\{ {e-1+i\choose e-1},\;
                         \sum_{j=i}^sc_j{e-1+j-i\choose e-1} \right\}\:.
  \end{align*}
  If equality holds for every $i$, then $\calQ/\calI$ is called
  \emph{compressed}, see \seccite[3]{RFrDLk84}.

  The number $n=\sum_{j=1}^sc_j$ is the type of $\calQ/\calI$ and,
  evidently, one has
  \begin{equation}
    \label{eq:cl}
    \begin{aligned}
      h_{\calQ/\calI}(i) &\dle \min\set{h_\calQ(i),nh_\calQ(s-i)}\\
      &\deq \min\left\{ {e-1+i\choose e-1},\; n{e-1+s-i\choose e-1}
      \right\}\:.
    \end{aligned}
  \end{equation}
\end{ipg}

\begin{prp}
  \label{prp:m}
  Adopt the notation from {\rm \pgref{huineq}} and let $e=3$. If the
  ideal $\calI$ has initial degree $d_{1,1}\ge 2$ and resolution
  format $(1,m,m+1,2)$ for $4\le m \le 7$, then one has
  $d_{3,2} > d_{1,1} + d_{1,2}$.
\end{prp}

As one has $d_{1,1} + d_{1,2} \ge 2d_{1,1}$ it follows that ideals
$\calI$ of format $(1,5,6,2)$, $(1,6,7,2)$, and $(1,7,8,2)$ avoid the
numerical obstruction to being licci from \pgref{huineq}.

\begin{prf*}
  To simplify the notation from \pgref{huineq}, set
  \begin{equation*}
    R=\calQ/\calI\,,\quad d=d_{1,1}\,, \quad d'=d_{1,2}\,, \quad\text{and}\quad s=d_{3,2}-3\:.
  \end{equation*}
  Notice that $R$ is Artinian of socle degree $s$.  The inequalities
  $2\le d_{1,1}$ and $d_{1,1}+2\le d_{3,2}$ yield
  \begin{equation*}
    1\le s \qand d \le s+1\:.
  \end{equation*}
  We assume that $s+3 \le d + d'$ holds and aim for a
  contradiction. To this end it is by \eqref{cl} sufficient to prove
  that there is a $u\le s$ such that the inequality
  $h_R(u) > 2h_{\calQ}(s-u)$ holds.

  \emph{Case 1.} Assume that $d < \frac{s+3}{2}$ holds. One then has
  $d' \ge s+3-d > d$, and hence
  \begin{align*}
    h_R(i) &\ge h_{\calQ}(i) - h_{\calQ}(i-d) - (m-1)h_{\calQ}(i-d')\\
           &\ge h_{\calQ}(i) - h_{\calQ}(i-d) - (m-1)h_{\calQ}(i-s-3+d)\\
           & \ge h_{\calQ}(i) - h_{\calQ}(i-d) - 6h_{\calQ}(i-s-3+d)\:.
  \end{align*}
  For $d=2$ one has
  \begin{equation*}
    h_R(s) \ge h_\calQ(s) - h_\calQ(s-2) = \binom{s+2}{2} - \binom{s}{2} = 2s+1 > 2 =
    2h_\calQ(0)\:.
  \end{equation*}
  Now let $d\ge 3$ and set $u = s + 2 - \lceil\frac{d}{2}\rceil$;
  notice that $u \le s$. Straightforward computations yield
  \begin{equation*}
    2h_{\calQ}(s-u) = 
    \begin{cases}
      \frac{1}{4}d(d-2) & \text{for $d$ even}\\
      \frac{1}{4}(d+1)(d-1) & \text{for $d$ odd}
    \end{cases}
  \end{equation*}
  and
  \begin{equation*}
    h_{\calQ}(u) - h_{\calQ}(u-d) - 6h_{\calQ}(u-s-3+d) =
    \begin{cases}
      \frac{1}{4}d(4s-7d+8) & \text{for $d$ even}\\
      \frac{1}{4}(d(4s-7d+12)+3) & \text{for $d$ odd.}
    \end{cases}
  \end{equation*}
  For $d$ even one now has
  \begin{align*}
    h_R(u) - 2h_{\calQ}(s-u) &\ge \frac{d(4s-7d+8)}{4} - \frac{d(d-2)}{4}\\
                             & = \frac{d(2(s-2d+3)-1)}{2}\\
                             & \ge \frac{d}{2}\:,
  \end{align*}
  and for $d$ odd
  \begin{align*}
    h_R(u) - 2h_{\calQ}(s-u) &\ge  \frac{d(4s-7d+12)+3}{4} - \frac{(d+1)(d-1)}{4}\\
                             & = d(s-2d+3) +1\\
                             & \ge d + 1\:.
  \end{align*}
  That is, independent of the parity of $d$ one has
  $h_R(u) > 2h_\calQ(s-u)$.

  \emph{Case 2.} Assume that $d \ge \frac{s+3}{2}$ holds; one has
  \begin{equation*}
    h_R(i) \ge h_{\calQ}(i) - mh_{\calQ}(i-d) \ge h_{\calQ}(i) - 7h_{\calQ}(i-d)\:.
  \end{equation*}
  Set $u = \lfloor\frac{s+d}{2}\rfloor$ and $\sigma = s-d+3$. Notice
  that one has
  \begin{equation*}
    1 \le \sigma \le d \le s+1 \qand u \le s\:.
  \end{equation*}
  Straightforward computations yield
  \begin{equation*}
    2h_{\calQ}(s-u) = 
    \begin{cases}
      \frac{1}{4}(\sigma^2-1) & \text{for $s+d$ even }\\
      \frac{1}{4}\sigma(\sigma+2) & \text{for $s+d$ odd}
    \end{cases}
  \end{equation*}
  and
  \begin{equation*}
    h_{\calQ}(u) - 7h_{\calQ}(u-d) =
    \begin{cases}
      \frac{1}{4}(2d(d+\sigma) - 3(\sigma^2-1))
      & \text{for $s+d$ even}\\
      \frac{1}{4}(2d(d+\sigma -1) - 3\sigma(\sigma-2)) & \text{for
        $s+d$ odd.}
    \end{cases}
  \end{equation*}

  The inequality $d \ge \sigma$ explains the second inequality in each
  of the computation below. For $s+d$ even one has
  \begin{align*}
    h_R(u) - 2h_{\calQ}(s-u) &\ge \frac{2d(d+\sigma) -
                               3(\sigma^2-1)}{4}
                               - \frac{\sigma^2-1}{4}\\
                             & \ge \frac{4\sigma^2 - 4(\sigma^2 -1)}{4}\\
                             & = 1\:,
  \end{align*}
  and for $s+d$ odd
  \begin{align*}
    h_R(u) - 2h_{\calQ}(s-u) &\ge \frac{2d(d+\sigma - 1)
                               -3\sigma(\sigma -2)}{4}
                               - \frac{\sigma(\sigma+2)}{4}\\
                             & \ge \frac{2\sigma(2\sigma - 1) - 4\sigma(\sigma-1)}{4}\\
                             & = \frac{\sigma}{2}\:.
  \end{align*}
  That is, independent of the parity of $d$ one has
  $h_R(u) > 2h_\calQ(s-u)$.
\end{prf*}

\begin{prp}
  \label{prp:n}
  Adopt the notation from {\rm \pgref{huineq}} and let $e=3$. If the
  ideal $\calI$ has initial degree $d_{1,1}\ge 2$ and resolution
  format $(1,5,n+4,n)$ for $1\le n \le 4$, then one has
  $d_{3,2} > d_{1,1} + d_{1,2}$.
\end{prp}

As one has $d_{1,1} + d_{1,2} \ge 2d_{1,1}$ it follows, in particular,
that ideals $\calI$ of format $(1,5,7,3)$ and $(1,5,8,4)$ avoid the
numerical obstruction to being licci from \pgref{huineq}.

\begin{prf*}
  As in the proof of \prpref{m}, set $R=\calQ/\calI$, $d=d_{1,1}$,
  $d'=d_{1,2}$, and $s=d_{3,2}-3$. Note that $R$ is Artinian of socle
  degree $s$, and that one has $1\le s$ and $d \le s+1$.  We assume
  that $s+3 \le d + d'$ holds and aim for a contradiction. To this end
  it is by \eqref{cl} sufficient to prove that there is a $u\le s$
  such that the inequality $h_R(u) > nh_{\calQ}(s-u)$ holds.

  \emph{Case 1.} Assume that $d \le \frac{s+3}{2}$ holds. One then has
  $d' \ge s+3-d \ge d$, and hence
  \begin{align*}
    h_R(i) &\ge h_{\calQ}(i) - h_{\calQ}(i-d) - 4h_{\calQ}(i-d')\\
           &\ge h_{\calQ}(i) - h_{\calQ}(i-d) - 4h_{\calQ}(i-s-3+d)\:.
  \end{align*}
  For $d=2$ one has $s \ge 2$ and
  $h_R(s) \ge h_\calQ(s) - h_\calQ(s-2) = 2s+1 > 4 \ge nh_\calQ(0)$ as
  computed in the proof of \prpref{m}.

  Now let $d\ge 3$ and set $u = s + 2 - \lceil\frac{d}{2}\rceil$;
  notice that $u \le s$. Straightforward computations yield
  \begin{equation*}
    nh_{\calQ}(s-u) \le 4h_{\calQ}(s-u) = 
    \begin{cases}
      \frac{1}{2}d(d-2) & \text{for $d$ even}\\
      \frac{1}{2}(d+1)(d-1) & \text{for $d$ odd}
    \end{cases}
  \end{equation*}
  and
  \begin{equation*}
    h_{\calQ}(u) - h_{\calQ}(u-d) - 4h_{\calQ}(u-s-3+d) =
    \begin{cases}
      \frac{1}{2}d(2s-3d+5) & \text{for $d$ even}\\
      \frac{1}{2}(d(2s-3d+6) + 1) & \text{for $d$ odd.}
    \end{cases}
  \end{equation*}
  For $d$ even one now has
  \begin{align*}
    h_R(u) - nh_{\calQ}(s-u) &\ge \frac{d(2s-3d+5)}{2} - \frac{d(d-2)}{2}\\
                             & = \frac{d(2(s-2d+3)+1)}{2}\\
                             & \ge \frac{3d}{2}\:,
  \end{align*}
  and for $d$ odd
  \begin{align*}
    h_R(u) - nh_{\calQ}(s-u) &\ge \frac{d(2s-3d+6) + 1}{2} - \frac{(d+1)(d-1)}{2} \\
                             & = \frac{2d(s-2d+3) + 2}{2}\\
                             & \ge d+1\:.
  \end{align*}
  That is, independent of the parity of $d$ one has
  $h_R(u) > nh_\calQ(s-u)$.

  \emph{Case 2.} Assume that $d > \frac{s+3}{2}$ holds; one has
  \begin{equation*}
    h_R(i) \ge h_{\calQ}(i) - 5h_{\calQ}(i-d)\:.
  \end{equation*}
  Set $u = \lfloor\frac{s+d}{2}\rfloor$ and $\sigma = s-d+3$. Notice
  that one has
  \begin{equation*}
    2 \le \sigma +1 \le d \le s+1 \qand u \le s\:.
  \end{equation*}
  Straightforward computations yield
  \begin{equation*}
    nh_{\calQ}(s-u) \le 4h_{\calQ}(s-u) = 
    \begin{cases}
      \frac{1}{2}(\sigma^2-1) & \text{for $s+d$ even }\\
      \frac{1}{2}\sigma(\sigma+2) & \text{for $s+d$ odd}
    \end{cases}
  \end{equation*}
  and
  \begin{equation*}
    h_{\calQ}(u) - 5h_{\calQ}(u-d) =
    \begin{cases}
      \frac{1}{2}(d(d+\sigma) - (\sigma^2-1))
      & \text{for $s+d$ even}\\
      \frac{1}{2}(d(d+\sigma -1) - \sigma(\sigma-2)) & \text{for $s+d$
        odd.}
    \end{cases}
  \end{equation*}

  The inequality $d \ge \sigma+1$ explains the second inequality in
  each of the computation below. For $s+d$ even one has
  \begin{align*}
    h_R(u) - nh_{\calQ}(s-u) &\ge \frac{d(d+\sigma) - (\sigma^2-1)}{2}
                               - \frac{\sigma^2-1}{2}\\
                             & \ge \frac{(\sigma+1)(2\sigma+1) - 2(\sigma^2 -1)}{2}\\
                             & = \frac{3(\sigma+1)}{2}\:,
  \end{align*}
  and for $s+d$ odd
  \begin{align*}
    h_R(u) - nh_{\calQ}(s-u) &\ge \frac{d(d+\sigma - 1) -
                               \sigma(\sigma -2)}{2}
                               - \frac{\sigma(\sigma+2)}{2}\\
                             & \ge \frac{(\sigma+1)2\sigma - 2(\sigma^2-\sigma)}{2}\\
                             & = 2\sigma\:.
  \end{align*}
  That is, independent of the parity of $d$ one has
  $h_R(u) > nh_\calQ(s-u)$.
\end{prf*}

\appendix
\section{Realizability of resolution formats}

\noindent
In the power series ring $\kk[\mspace{-2mu}[X,Y,Z]\mspace{-2mu}]$
there exists by \prpcite[6.2]{DABDEs77} an ideal of resolution format
$(1,m,m,1)$, i.e.\ a Gorenstein ideal of grade $3$, for every odd
$m \ge 3$. We start by recording that this is the case in any local
ring of sufficient depth.

\begin{lem}
  \label{lem:gor}
  Let $Q$ be a local ring of depth at least $3$. For every odd
  $m \ge 3$ there exists a grade $3$ perfect ideal of resolution
  format $(1,m,m,1)$.
\end{lem}

\begin{prf*}
  Denote by $\fn$ the maximal ideal of $Q$ and let $x,y,z$ be a
  regular sequence in $\fn$.  Let $l$ be an integer with $m=2l+1$ and
  define a skew symmetric matrix $V_l$ as in \seccite[3]{CVW-1}. Let
  $I$ be the ideal generated by the sub-maximal Pfaffians of $V_l$. It
  follows from Proposition 3.3 in \emph{loc.\,cit.} that the radical
  of $I$ contains $x$, $y$, and $z$, so $I$ has grade (at least) $3$,
  whence by \thmcite[2.1]{DABDEs77} it is a Gorenstein ideal of format
  $(1,m,m,1)$.
\end{prf*}

\begin{thm}
  Let $(Q,\fn)$ be a local ring of depth at least $3$. For every
  resolution format $\ff = (1,m,m+n-1,n)$ with $m\ge 3$ and $n\ge 1$
  that is not $(1,m,m,1)$ with $m$ even and not $(1,3,n+2,n)$ with
  $n\ge 2$ there exists a grade $3$ perfect ideal $I$ in $Q$ of
  resolution format $\ff$.
\end{thm}

\begin{prf*}
  By \lemref{gor} there exists for every odd $g\ge 3$ a grade $3$
  perfect ideal with resolution format $(1,g,g,1)$. Fix a format
  $\ff = (1,m,m+n-1,n)$ with $m\ge 3$ and $n\ge 2$. If $m=4$ and $n$
  is even, then it follows from \prpref{link}(b) that there exists an
  ideal of format $\ff$ that is linked to one of format
  $(1,n+3,n+3,1)$. If $m=4$ and $n$ is odd, then it follows from
  \prpref{link}(a) that there exists an ideal of format $\ff$ that is
  linked to one of format $(1,n+2,n+2,1)$. If $m\ge 5$ is odd, then it
  follows from $n-1$ applications of \corref{link}(b) that there is an
  ideal of format $\ff$ that is linked to one of format
  $(1,m,m,1)$. If $m\ge 6$ is even, then it follows from
  \corref{link}(a) and $n-2$ applications of \corref{link}(b) that
  there exists an ideal of format $\ff$ that is linked to one of
  format $(1,m-1,m-1,1)$.
\end{prf*}

\def\soft#1{\leavevmode\setbox0=\hbox{h}\dimen7=\ht0\advance \dimen7
  by-1ex\relax\if t#1\relax\rlap{\raise.6\dimen7
    \hbox{\kern.3ex\char'47}}#1\relax\else\if T#1\relax
  \rlap{\raise.5\dimen7\hbox{\kern1.3ex\char'47}}#1\relax \else\if
  d#1\relax\rlap{\raise.5\dimen7\hbox{\kern.9ex
      \char'47}}#1\relax\else\if D#1\relax\rlap{\raise.5\dimen7
    \hbox{\kern1.4ex\char'47}}#1\relax\else\if l#1\relax
  \rlap{\raise.5\dimen7\hbox{\kern.4ex\char'47}}#1\relax \else\if
  L#1\relax\rlap{\raise.5\dimen7\hbox{\kern.7ex
      \char'47}}#1\relax\else\message{accent \string\soft \space #1
    not defined!}#1\relax\fi\fi\fi\fi\fi\fi}
\providecommand{\MR}[1]{\mbox{\href{http://www.ams.org/mathscinet-getitem?mr=#1}{#1}}}
\renewcommand{\MR}[1]{\mbox{\href{http://www.ams.org/mathscinet-getitem?mr=#1}{#1}}}
\providecommand{\arxiv}[2][AC]{\mbox{\href{http://arxiv.org/abs/#2}{\sf
      arXiv:#2 [math.#1]}}} \def\cprime{$'$}
\providecommand{\bysame}{\leavevmode\hbox to3em{\hrulefill}\thinspace}
\providecommand{\MR}{\relax\ifhmode\unskip\space\fi MR }
\providecommand{\MRhref}[2]{%
  \href{http://www.ams.org/mathscinet-getitem?mr=#1}{#2} }
\providecommand{\href}[2]{#2}

\end{document}